\documentclass[numreferences]{kluwer}

\def\proofname{Proof}
\newproof{pf}{\proofname}

\def\R{{\mathbb R}}

\def\sign{\mathop{\hbox{sign}}}
\def\normo#1{\left\|#1\right\|}
\def\modo#1{\left|#1\right|}
\def\Ave{\mathop{\hbox{\rm Ave}}\limits}

\newtheorem{Theorem}{Theorem}
\newtheorem{Conjecture}[Theorem]{Conjecture}
\newtheorem{Lemma}[Theorem]{Lemma}

\begin{document}
\begin{article}
\begin{opening}

\title{Embeddings of rearrangement invariant spaces
that are not strictly singular}

\author{S.J.\surname{Montgomery-Smith}\email{stephen@math.missouri.edu}
\thanks{Research supported in part by grants from the N.S.F.\ and the
Research Board of the University of Missouri.}}
\institute{University of Missouri}
\author{E.M.\surname{Semenov}\email{root@func.vsu.su}
\thanks{Research supported in part by a grant from the N.S.F.\
and a grant from Russia.}}
\institute{Voronezh State University}

\runningtitle{Embeddings of r.i. spaces}
\runningauthor{Montgomery-Smith and Semenov}

\begin{ao}
Department of Mathematics,
University of Missouri, 
Columbia, MO 65211,
U.S.A.
or 
Department of Mathematics,
Voronezh State University, Universitetskaya pl.1,
Voronezh 394693, Russia.
\end{ao}

\begin{abstract}
We give partial answers to the following conjecture:
the natural embedding of a rearrangement invariant
space $E$ into $L_1([0,1])$ is strictly singular if and only if
$G$ does not
embed into $E$ continuously, where $G$ is the closure of the simple functions
in the
Orlicz space $L_\Phi$ with $\Phi(x) = \exp(x^2)-1$.
\end{abstract}


{Primary 46E30, 47B38; Secondary 60G50}

\end{opening}

In this paper we ask the following question.  Given a rearrangement
invariant space $E$ on $[0,1]$, when is the natural embedding
$E \subset L_1([0,1])$ strictly singular.  
(We refer the reader to \cite{LT} for the definition and properties
of rearrangement invariant spaces.)
We define a linear map between
two normed spaces to be {\it strictly singular\/} if there does not
exist an infinite dimensional subspace of the domain upon which
the operator is an isomorphism.

This question is a natural extension of similar work by del Amo,
Hern\'andez, S\'anchez and Semenov \cite{AH}, 
when they considered the problem
of which embeddings between rearrangement invariant spaces are not
disjointly strictly singular.  A positive linear operator between two Banach
lattices is {\it disjointly strictly singular\/}
if there exists an
infinite sequence of non-zero disjoint elements in the domain
such that the operator is an isomorphism on the span of this sequence.
This work \cite{AH} contains a number of very sharp results, giving some very
clear criteria.

However the question concerning when such maps are strictly singular seems
to be more difficult.  For this reason, we will restrict ourselves to
considering the case when the range is $L_1([0,1])$.  Even then, we do
not have complete answers, and in this paper, we leave as many questions
unanswered as we answer.

The ``other end'' was investigated by Novikov \cite{N}, who showed that
the natural embedding $L_\infty([0,1]) \subset E$ is strictly singular
unless $E$ is equivalent to $L_\infty([0,1])$.  (The case when 
$E = L_p([0,1])$ follows from a classical result of Grothendieck,
Theorem 5.2 of \cite{R}.)

In answering our question,
there is one rearrangement invariant space that plays a prominent role.
This rearrangement invariant
space, denoted by $G$, 
is the closure of the simple functions
in the Orlicz space corresponding to the Orlicz function $e^{x^2}-1$.  
(We define Orlicz spaces below.)
The reason why this space plays such a role is as follows.
Consider the Rademacher functions on $[0,1]$ given by
$r_n(t) = \sign(\sin(2^n \pi t))$.  The following
result of Rodin and Semenov \cite{RS} is well known.

\begin{Theorem}  
\label{rodin-sem}
Let $E$ be a rearrangement invariant space on $[0,1]$.
Then the following are equivalent.
\begin{enumerate}
\item The sequence
$(r_n)$ in $E$ is equivalent to the unit vector basis of $\ell_2$;
\item $G$ embeds continuously into $E$;
\item There is a constant $c>0$ such that
$\normo{\sum_{i=1}^n r_i}_E \le c \sqrt n$.
\end{enumerate}
\end{Theorem}

It is clear that 
if $E$ is a rearrangement invariant space on $[0,1]$
that contains $G$ continuously, then the natural embedding
$E \subset L_1$ is not strictly singular.  Here, of course, the subspace
on which the norms are equivalent is the span of the Rademacher functions.

For this reason, it is natural to pose the following conjecture.

\begin{Conjecture}  
\label{ss-conj}
Let $E$ be a rearrangement invariant space on $[0,1]$.
Then the natural embedding $E \subset L_1$ is not strictly
singular if and only if $G$ embeds into $E$ continuously.
\end{Conjecture}

It will become
apparent below that the following conjecture implies the previous
one.

\begin{Conjecture}
\label{sign-conj}  
Given $x_1,\dots,x_n \in L_1([0,1])$, 
and a rearrangement invariant space $E$ on $[0,1]$,
there exists signs $\epsilon_1,\dots,\epsilon_n = \pm 1$ such that
$$ \normo{\sum_{i=1}^n \epsilon_i x_i}_{L_1}
   \ge
   c^{-1}
   \normo{\sum_{i=1}^n r_i\normo{x_i}_1}_E .$$
\end{Conjecture}

Unfortunately we are not able to prove either of these conjectures
without some additional hypotheses.  

First,
let us introduce some examples of rearrangement invariant spaces.
A function $\Phi:\R\to[0,\infty)$ is called an {\it Orlicz function\/}
if it is convex, even, and takes zero to zero.  The {\it Orlicz space\/}
$L_\Phi$ is the collection of all equivalence classes of
measurable functions (where the equivalence relation is equal almost
everywhere) on $[0,1]$ such that
the norm:
$$ \normo x_\Phi = \inf\left\{ \lambda: \int_0^1 \Phi(x(t)/\lambda) \, dt
   \le 1 \right\} $$
is finite.

Another class of examples is the Lorentz spaces.  If $\varphi:[0,1]\to[0,1]$ is
increasing and concave with $\varphi(0)=0$ and $\varphi(1)=1$, then
the {\it Lorentz space\/} 
$\Lambda(\varphi)$ consists of all equivalence classes of
measurable functions on $[0,1]$ for
which the norm
$$ \normo x_{\Lambda(\varphi)}
   =
   \int_0^1 x^*(t) \, d\varphi(t) $$
is finite.  Here, as in the rest of the paper, $x^*$ denotes the non-increasing
rearrangement of $\modo x$.

One more class of examples is the Marcinkiewicz spaces.  
If $\varphi:[0,1]\to[0,1]$ is
as above,
then the {\it Marcinkiewicz space\/} $M(\varphi) = \Lambda(\varphi)^*$ 
consists of all equivalence classes of measurable functions on
$[0,1]$ for which the norm
$$ \normo x_{M(\varphi)}
   =
   \sup_{0<t<1} {\int_0^t x^*(s) \, ds \over \varphi(t)} $$
is finite.

By definition, $G$ 
is the closure of the simple functions
in the Orlicz space corresponding to the Orlicz function $e^{x^2}-1$.  
However, it can also be shown that it has an 
equivalent Marcinkiewicz norm, that is, there is a constant $c>0$ such that
$c^{-1} \normo x_G \le \normo x_{M(\varphi)} \le c \normo x_G$, where
$\varphi(t) = \displaystyle {t \over \sqrt{\log(e/t)}}$.

Another space that will concern us is the space $G_1 = \Lambda(\varphi)$
where $\varphi(t) = \displaystyle {2 \over \sqrt{\log(e^2/t)}}$.  It is a
simple matter to show that $G_1$ embeds continuously into $G$.  
Furthermore, for characteristic functions,
the norms on $G$ and $G_1$ coincide.

Let us now consider the
concept of $D$-convex rearrangement invariant spaces.  This notion was
introduced by Kalton \cite{K}, and studied extensively by Montgomery-Smith
and Semenov \cite{MSe}, where several equivalent properties were given.  Perhaps
the easiest definition to work with is the following.  We will say that
a rearrangement invariant space $E$ is {\it $D$-convex\/} if there is a
family of Orlicz functions $\Phi_\alpha:\R \to [0,\infty)$, and a constant
$c>0$, such that
$$ c^{-1} \normo x_E \le  \sup_\alpha \normo x_{\Phi_\alpha} \le \normo x_E .$$

Note that the Marcinkiewicz spaces are
$D$-convex, because for each $0<t<1$ the map
$x \mapsto \int_0^t x^*(s) \, ds$ is equivalent to the Orlicz norm
given by the Orlicz function $\Phi(s) = (s-t^{-1})^+$.

The best result we have obtained so far regarding 
Conjectures~\ref{ss-conj} and~\ref{sign-conj} is the following.

\begin{Theorem}
\label{ss-sign}
Conjectures~\ref{ss-conj} and~\ref{sign-conj} are true if $E$ is $D$-convex.
\end{Theorem}

From here we are able to get a weaker version of Conjecture~\ref{ss-conj}.

\begin{Theorem}  
\label{g1}
Let $E$ be a rearrangement invariant space on $[0,1]$.
If the natural embedding $E \subset L_1$ is not strictly
singular, then $G_1$ embeds into $E$ continuously.
\end{Theorem}

We proceed with the proofs.

\begin{pf*}{Proof of first part of Theorem~\ref{ss-sign}}
It is sufficient to consider the case when $E$ is an Orlicz space
$L_\Phi$, where $\Phi$ is an Orlicz function.  Suppose that
$$ \sup_\epsilon \normo{\sum_{i=1}^n \epsilon_i x_i}_\Phi \le 1 ,$$
where we write $\epsilon = (\epsilon_1,\dots,\epsilon_n)$.
Thus
$$ \sup_\epsilon \int \Phi\left(\sum_{i=1}^n \epsilon_i x_i(s)
   \right) \, ds \le 1 .$$
Let
$$ F_\epsilon = \{s\in[0,1] : \sign(x_i(s)) = \epsilon_i\} .$$
If $\epsilon = (\epsilon_1,\dots,\epsilon_n)$
and
$\eta = (\eta_1,\dots,\eta_n)$, let us set
$\eta \epsilon = (\eta_1 \epsilon_1,\dots,\eta_n \epsilon_n)$.
Then
\begin{eqnarray*}
   \Ave_\eta \Phi\left(\sum_{i=1}^n \eta_i \normo{x_i}_1\right)
   &=&
   \Ave_\eta \Phi\left(\sum_\epsilon \int I_{F_{\eta \epsilon}}(s)
   \sum_i \eta_i \modo{x_i(s)} \, ds \right) \\
   &=&
   \Ave_\eta \Phi\left(\sum_\epsilon \int I_{F_{\eta \epsilon}}(s)
   \sum_i \epsilon_i x_i(s) \, ds \right) \\
   &\le&
   \Ave_\eta \int \Phi\left(\sum_\epsilon I_{F_{\eta \epsilon}}(s)
   \sum_i \epsilon_i x_i(s) \right) \, ds \\
   &=&
   \Ave_\eta \int \sum_\epsilon I_{F_{\eta \epsilon}}(s)
   \Phi\left(
   \sum_i \epsilon_i x_i(s) \right) \, ds \\
   &=&
   \Ave_\epsilon \sum_\eta \int I_{F_{\eta \epsilon}}(s)
   \Phi\left(
   \sum_i \epsilon_i x_i(s) \right) \, ds \\
   &\le&
   \sup_\epsilon \int
   \Phi\left(
   \sum_i \epsilon_i x_i(s) \right) \, ds \\
   &\le&
   1 . 
\end{eqnarray*}
Therefore
$$ \normo{\sum_i r_i \normo{x_i}_1}_\Phi \le 1 .$$
\end{pf*}\qed

\begin{pf*}{Proof of second part of Theorem~\ref{ss-sign}}
Suppose that there is an infinite dimensional subspace $F \subset E$ such that
the norms of $L_1$ and $E$ are equivalent on $F$.  By Dvoretzky's Theorem 
(see \cite{MSc} Chapter~4),
for each integer $n$, there is an $n$-dimensional subspace $H$ of $F$ such
that $H$ is $2$-isomorphic to Hilbert space.  Pick
an orthonormal basis $x_1,\dots,x_n$ for $H$.  Then there exist signs
$\epsilon_1,\dots,\epsilon_n = \pm 1$ such that
$$ \normo{\sum_i r_i}_E
   \le c \,
   \normo{\sum_i \epsilon_i x_i}_H
   = c \sqrt{n} .$$
Thus the result follows by Theorem~\ref{rodin-sem}.
\end{pf*}\qed

In order to prove Theorem~\ref{g1}, we need the following Lemma.  This result
may also be found in \cite{KPS} (Theorem~2.5.7.).

\begin{Lemma}
\label{M-embed}
Given a rearrangement invariant space $E$, there exists
an increasing function $\varphi:[0,1]\to[0,1]$ such that $\normo x_E \ge
\normo x_{M(\varphi)}$, but if $x$ takes only values $0$ or $1$, then
$\normo x_E = \normo x_{M(\varphi)}$.
\end{Lemma}

\begin{pf}
Let 
$$ \varphi(t) = {t\over \normo{I_{[0,t]}}_E} .$$
The latter property is obvious.  To show the former, recall (see for example
\cite{LT}) the
space $E'$ to be those functions $y$ on $[0,1]$ for which the
norm
$$ \normo y_{E'} = \sup\left\{ \int x(t) y(t) \, dt : \normo x_E \le 1
   \right\}$$
is finite.  Notice that
$$ \normo{I_{[0,t]}}_E \cdot \normo{I_{[0,t]}}_{E'} = t ,$$
that is, $\normo{I_{[0,1]}}_{E'} = \varphi(t) $.  Then the result follows
since
$$ \int_0^t x^*(s) \, ds
   \le \normo x_E \normo{I_{[0,1]}}_{E'} .$$
\end{pf}\qed

\begin{pf*}{Proof of Theorem~\ref{g1}}
Suppose that the embedding $E \subset L_1([0,1])$ is not strictly singular.
Produce $\varphi$ as in Lemma~\ref{M-embed}.  Then it is clear that the embedding
$M(\varphi) \subset L_1([0,1])$ is also not strictly singular.  Hence
by Theorem~\ref{ss-sign}, $G$ embeds continuously into $M(\varphi)$.

In particular, this means that there is a constant $c>0$ such that for
all $t \in [0,1]$ we have
$$ \normo{I_{[0,t]}}_E = \normo{I_{[0,t]}}_{M(\varphi)}
   \le c \normo{I_{[0,t]}}_G
   \le c \psi(t) ,$$
where $\psi(t) = \displaystyle {2 \over \sqrt{\log(e^4/t)}} $.
Now, writing $x^* = \int_0^1 I_{[0,t]} \, d(-x^*(t))$, we obtain that
\begin{eqnarray*}
   \normo x_E
   &\le&
   \int_0^1 \normo{I_{[0,t]}}_E \, d(-x^*(t)) \\
   &\le& c \int_0^1 \psi(t) d(-x^*(t))
   = c \int_0^1 x^*(t) d\psi(t) 
   = c \normo x_{G_1} .
\end{eqnarray*}
\end{pf*}\qed

\end{article}
\end{document}